\documentclass[oneside,draft,12pt]{amsart}
\def\date #1{\rightline{\vbox{\hsize 1.5true in\noindent #1}}\bigskip}

\def\block #1{\vbox{\hsize 2.5 true cm\noindent \bf#1:}}

\def\memo TO:#1FROM:#2SUBJECT:#3DATE:#4\par{\centerline{\sl MEMO!}
      \bigskip \hrule height1pt \medskip
      \vbox{\parindent=75pt\parskip=1pt
      \item{\block{TO}}#1
      \item{\block{FROM}}#2
      \item{\block{SUBJECT}}#3
      \item{\block{DATE}}#4 }\medskip \hrule height1pt \bigskip}

\newcount\spk
\def\beginscript {\bgroup \parindent=0pt \spk=1 \sl \rightskip.4in
      \def\par {\ifnum\spk=1 \endgraf \it \spk=2 \leftskip.4in \rightskip0in
               \else \endgraf \sl \spk=1 \leftskip0in \rightskip.4in \fi}}

\def\endscript {\egroup}

\newcount\var \newcount\pw \newcount\tmp \newcount\cnt
\def\pow#1#2#3{\var=#1 \pw=#2 \tmp=\var \cnt=1
      \loop \multiply\var by\tmp \advance\cnt by1 \ifnum\cnt<\pw \repeat
      \global#3=\var}

\newcount\xf \newcount\xnf \newcount\yf \newcount\ynf \newcount\zf \newcount\znf
\def\n {\number}

\def\fermat#1#2#3{$\global\xf=#1 \global\yf=#2 \global\pw=#3
      \pow{\xf}{\pw}{\xnf} \pow{\yf}{\pw}{\ynf}
      \global\tmp=\ynf \global\advance\tmp by\xnf
      {\n\xf}^{\n\pw}+{\n\yf}^{\n\pw}={\n\tmp}$.\hfil\break
      \ifnum\xf>\yf \zf=\xf \else \zf=\yf \fi
      \loop {\pow{\zf}{\pw}{\znf}} \ifnum\znf<\tmp \advance\zf by1 \repeat
      \ifnum\znf=\tmp  The sum seems to be exactly ${\n\zf}^{\n\pw}$.
        \ifnum\pw=2 {\it Yawn!\/} Tell me something I don't know, will you?
        \else Incredible!  But, perhaps you'd better check my work.\fi
      \else \advance\zf by-1
      This lies between ${\n\zf}^{\n\pw}$\pow{\zf}{\pw}{\znf}($={\n\znf}$)
      \advance\zf by1
      and ${\n\zf}^{\n\pw}$\pow{\zf}{\pw}{\znf}($={\n\znf}$).\fi}%

\def\frame #1#2#3#4{\vbox{\hrule height #1pt%    TOP RULE
 \hbox{\vrule width #1pt\kern #2pt%              RULE AND SPACE ON LEFT
 \vbox{\kern #2pt%                               SPACE AT TOP
 \vbox{\hsize #3\noindent #4}%                   MATERIAL THAT WILL BE BOXED
 \kern #2pt}%                                    SPACE AT BOTTOM
 \kern #2pt\vrule width #1pt}%                   SPACE AND RULE ON RIGHT
 \hrule height0pt depth #1pt}}%,                 BOTTOM RULE

\def\fitframe #1#2#3{\vbox{\hrule height#1pt%    TOP RULE
 \hbox{\vrule width#1pt\kern #2pt%               RULE AND SPACE ON LEFT
 \vbox{\kern #2pt\hbox{#3}\kern #2pt}%           TOP, MATERIAL, BOTTOM
 \kern #2pt\vrule width#1pt}%                    SPACE AND RULE ON RIGHT
 \hrule height0pt depth#1pt}}%                   BOTTOM RULE

\def\shframe #1#2#3#4{\vbox{\hrule height 0pt%   NO TOP SHADOW
 \hbox{\vrule width #1pt\kern 0pt%               LEFT SHADOW
 \vbox{\kern-#1pt\frame{.3}{#2}{#3}{#4}%         SHADOW STARTS #1 PT FROM TOP
 \kern-.3pt}%                                    MOVE UP RULE THICKNESS AT BOT.
 \kern-#1pt\vrule width 0pt}%                    STOPS #1 PT FROM RT; NO RT SHAD
 \hrule height #1pt}}%                           BOTTOM SHADOW

\newcount\tw        %top width in pt
\newcount\bw        %bottom width in pt
\newcount\h         %height in pt
\newcount\bs        %bottom shift in pt
\newcount\th        %line thickness in 1/64th of a pt
\newcount\gp        %line gap in 1/64th of a point
\newcount\rs        %running shift
\newcount\rw        %running width
\newcount\rh        %running height
\newcount\tmp       %for temporarily storing variables

\def\trap #1#2#3#4#5#6{\vbox{\offinterlineskip
      \tw=#1 \bw=#2 \h=#3 \bs=#4 \th=#5 \gp=#6 \rh=0
      \multiply\tw by 65536 \multiply\bw by 65536 \multiply\bs by 65536
      \multiply\th by 1024 \multiply\gp by 1024
      \loop
      \tmp=\bs \multiply\tmp by\rh \divide\tmp by\h
      \rs=\tmp                                        %running shift calculated
      \tmp=\bw \advance\tmp by-\tw \multiply\tmp by\rh
      \divide\tmp by\h \advance\tmp by\tw \rw=\tmp    %running width calculated
      \hbox{\kern\rs sp\vrule height0sp depth\th sp width\rw sp}%RULE DRAWN HERE
      \vskip\gp sp                                              %GAP LEFT HERE
      \ifnum\rh<\h  \tmp=\rh \multiply\tmp by 65536
       \advance\tmp by\th \advance\tmp by\gp \divide\tmp by65536 \rh=\tmp
      \repeat}}%

\newcount\colnumber \newbox\col \newdimen\tmpdim \newdimen\size
\newdimen\coljump \coljump=.2 true in                  %GAP BETWEEN COLUMNS

\def\nstrut {\vrule height\topskip depth0pt width0pt}  %TO PROP THINGS UP
\def\divider{\hskip\coljump}
\def\dividerule{\dimen0=.4pt \dimen1=\coljump \advance\dimen1 by-\dimen0
      \divide\dimen1 by2
      \def\divider{\hskip\dimen1 \vrule width\dimen0 \hskip\dimen1}}

\def\beginart #1/#2{\vbox\bgroup#1 \colnumber=#2  \parskip=0pt
      \advance\colnumber by-1 \tmpdim=\coljump \multiply\tmpdim by\colnumber
      \size=\hsize  \advance\size by-\tmpdim
      \advance\colnumber by1 \divide\size by\colnumber
      \vbadness=10000 \hbadness=2000 \tolerance=2000
      \setbox\col=\vbox\bgroup\hsize\size \noindent\nstrut}

\def\endart {\global\size=\baselineskip \vfil \egroup
      \multiply\size by\colnumber  \advance\size by-\baselineskip
      \tmpdim=\ht\col  \advance\tmpdim by\size  \divide\tmpdim by\colnumber
      \hbox{\splittopskip=\topskip \doittoit}\egroup}

\def\doittoit{\ifnum\colnumber>0 \vsplit\col to \tmpdim
      \global\advance\colnumber by-1
      \ifnum\colnumber>0 \divider \fi \doittoit \fi}

\def\chartable #1{\smallbreak\vbox{\noindent%
     \underbar{Characters in the {\sl#1\/} font, with decimal codes:}
     \raggedright \hbadness5000 \tolerance10000 \medskip
     \font\ft=#1 \ft \dimen0=14pt \baselineskip=\dimen0
     \ifdim\dimen0<3.25ex \baselineskip=3.25ex \fi \count255=0 \dimen0=10pt
     \loop \setbox0=\hbox{\char\count255}
     \ifdim\wd0>\dimen0 \dimen0=\wd0 \fi \advance\count255 by1
     \ifnum\count255<128 \repeat \count255=0
     \advance\dimen0 by 25pt \noindent \loop
     \hbox to\dimen0{\hbox to23pt{\hfil\rm\the\count255:\ }\char\count255\hfil}
     \advance\count255 by1 \ifnum\count255<128 \quad\repeat
     \smallbreak} \font\tenrm=cmr10 }

\def\ignore {\count255=0 \begingroup
      \loop \catcode\count255=14  % Make everything a comment character.
         \advance\count255 by1 \ifnum\count255<127
      \repeat \catcode`\!=0 }     % Makes ! an escape character.
{\catcode`\!=0 !gdef!E{!endgroup}}% Defines the `stop ignoring' command.

\def\stacksymbols #1#2#3#4{\def\theguybelow{#2}
    \def\verticalposition{\lower#3pt}
    \def\spacingwithinsymbol{\baselineskip0pt\lineskip#4pt}
    \mathrel{\mathpalette\intermediary#1}}
\def\intermediary#1#2{\verticalposition\vbox{\spacingwithinsymbol
      \everycr={}\tabskip0pt
      \halign{$\mathsurround0pt#1\hfil##\hfil$\crcr#2\crcr
               \theguybelow\crcr}}}

\def\monthname {\ifcase\month\or January\or February\or March\or
April\or May\or June\or July\or August\or September\or October\or
November\or December\fi}

\newcount\mins  \newcount\hours  \hours=\time \mins=\time
\def\now{\divide\hours by60 \multiply\hours by60 \advance\mins by-\hours
     \divide\hours by60         % NOTE: \divide only gives integer answers.
     \ifnum\hours>12 \advance\hours by-12
       \number\hours:\ifnum\mins<10 0\fi\number\mins\ P.M.\else
       \number\hours:\ifnum\mins<10 0\fi\number\mins\ A.M.\fi}

\newdimen\tempdim                 % For temporary storage.
\newdimen\othick   \othick=.4pt   % To set the outer rule thickness.
\newdimen\ithick   \ithick=.4pt   % To set the inner rule thickness.
\newdimen\spacing  \spacing=9pt   % To set the interline spacing.
\newdimen\abovehr  \abovehr=6pt   % Space above horizontal rules.
\newdimen\belowhr  \belowhr=8pt   % Space below horizontal rules.
\newdimen\nexttovr \nexttovr=8pt  % Space next to vertical rules.

\def\rr{\hfil\down{\abovehr}&\omit\vrsp\vrule width\othick\cr
     \noalign{\hrule height\ithick}\up{\belowhr}&}% To draw an \hrule.
\def\up#1{\tempdim=#1\advance\tempdim by1ex
     \vrule height\tempdim width0pt depth0pt}%   For space above a line.
\def\down#1{\vrule height0pt depth#1 width0pt}%  For space below a line.
\def\large#1#2{\setbox0=\vtop{\hsize#1 \lineskiplimit=0pt \lineskip=1pt
     \baselineskip\spacing \advance\baselineskip by 3pt \noindent
     #2}\tempdim=\dp0\advance\tempdim by\abovehr\box0\down{\tempdim}}
\def\vrsp{\hskip\nexttovr\relax}
\def\toprule#1{\def\startrule{\hrule height#1\relax}} % Set a top rule.
\toprule{\othick}                      % Picking the `\toprule' default.
\def\nstrut{\vrule height\spacing depth3.5pt width0pt}
\def\preamble#1{\def\startup{#1}}      % For `customized' preambles.
\preamble{&##}                         % Choosing the default preamble.
{\catcode`\!=\active
 \gdef!{\hfil\vrule width0pt\vrsp\vrule width\ithick\relax\vrsp&}}
\def\table #1{\vbox\bgroup \setbox0=\hbox{#1}
     \vbox\bgroup\offinterlineskip  \catcode`\!=\active
     \halign\bgroup##\vrule width\othick\vrsp&\span\startup\nstrut\cr
     \noalign{\medskip}
     \noalign{\startrule}\up{\belowhr}&}

\def\caption #1{\down{\abovehr}&\omit\vrsp\vrule width\othick\cr
     \noalign{\hrule height\othick}\egroup\egroup \setbox1=\lastbox
     \tempdim=\wd1 \hbox to\tempdim{\hfil \box0 \hfil} \box1 \smallskip
     \hbox to\tempdim{\advance\tempdim by-20pt\hfil\vbox{\hsize\tempdim
     \noindent #1}\hfil}\egroup}

\let\cc=\catcode
{\cc`\^^M=\active %
\gdef\losenolines{\cc`\^^M=\active \def^^M{\leavevmode\endgraf}}}
\def\literal {\begingroup \cc`\\=12 \cc`\{=12 \cc`\}=12 \cc`\$=12 \cc`\&=12
 \cc`\#=12 \cc`\%=12 \cc`\~=12 \cc`\_=12 \cc`\^=12 \cc`\*=12 \cc`\@=0
 \cc`\`=\active \losenolines \obeyspaces \tt}%
{\obeyspaces\gdef {\hglue.5em\relax}}

{\cc`\`=\active \gdef`{\relax\lq}}

\def\vquotingon{\cc`\"=\active}
\def\vquotingoff{\cc`\"=12}
\vquotingon
\def"{\literal\leavevmode\hbox\bgroup\com}
{\cc`\@=0 \cc`\\=12 @cc`@^^M=@active %
 @gdef@com#1"{#1@egroup@endgroup} %
 @gdef@thatisit^^M#1\endliteral{#1@endgroup@smallskip}}
\vquotingoff

\def\pattern #1#2{\count0=0
      \loop #1\advance\count0 by 1 \ifnum\count0<#2 \repeat}

\usepackage{amssymb, amscd}
\usepackage[all]{xy}

\theoremstyle{plain}
\newtheorem{theorem}{Theorem}[section]
\newtheorem{lemma}[theorem]{Lemma}

\newtheorem{proposition}[theorem]{Proposition}

\newtheorem{corollary}[theorem]{Corollary}

\theoremstyle{definition}

\newtheorem{definition}[theorem]{Definition}

\newtheorem{example}[theorem]{Example}

\theoremstyle{remark}
\newtheorem{remark}[theorem]{Remark}

\DeclareMathOperator{\holim}{holim}

\DeclareMathOperator{\aut}{Aut}

\newcommand{\bbF}{\ensuremath{\mathbb{F}}}

\newcommand{\bbR}{\ensuremath{\mathbb{R}}}

\newcommand{\bbC}{\ensuremath{\mathbb{C}}}
\newcommand{\bbZ}{\ensuremath{\mathbb{Z}}}
\newcommand{\bbS}{\ensuremath{\mathbb{S}}}

\newcommand{\fib}{\vrule width 1.75cm height .1pt}

\newcommand{\ra}{\ensuremath{\rightarrow}}

\renewcommand{\a}{\alpha}
\renewcommand{\b}{\beta}

\begin{document}

\title{Periodic Complexes and Group Actions}

\author{Alejandro Adem}

\address{Department of Mathematics\\
         University of Wisconsin\\
         Madison, Wisconsin 53706}

\email{adem@math.wisc.edu}

\thanks{Both authors were partially supported by grants from the
NSF and the CRM--Barcelona}

\author{Jeff H. Smith}
\address{Department of Mathematics\\
         Purdue University\\
         West Lafayette, Indiana 47907}
\email{jhs@math.purdue.edu}

\begin{abstract}
  In this paper we show that the cohomology of a connected
  $CW$--complex is periodic if and only if it is the base space of a
  spherical fibration with total space that is homotopically finite
  dimensional.  As applications we characterize those discrete groups
  that act freely and properly on $\mathbb R^n\times \mathbb \bbS^m$;
  we construct non--standard free actions of rank two simple groups on
  finite complexes $Y\simeq \bbS^n\times\bbS^m$; and we prove that a
  finite $p$--group $P$ acts freely on such a complex if and only if
  it does not contain a subgroup isomorphic to $(\bbZ/p)^3$.
\end{abstract}

\maketitle

\section{Introduction}

Let $G$ be a finite group for which all abelian subgroups
are cyclic. Swan \cite{Sw} has proved that any such group acts freely
on a finite complex $X\simeq\bbS^m$ for some $m>0$. These groups are
best described as those having {\sl periodic} group cohomology. 
Since then there have been numerous attempts to extend this type of result
to other classes of groups, including (1) infinite groups with periodic
cohomology and (2) finite groups with non--periodic cohomology.

In
this paper we prove a homotopy--theoretic characterization of
{\sl cohomological periodicity}, a notion which we now make precise.

\begin{definition}
A CW--complex $X$ has \emph{periodic} cohomology if there is a
cohomology class $\a\in H^*(X,\bbZ)$ with $|\a|>0$ and an
integer $d\ge 0$ such that the map
\[
\a\cup-\,\colon H^n_{loc}(X,\mathcal B)\to H^{n+|\a|}_{loc}(X,\mathcal
B)
\]
induced by the cup product with $\a$ is an isomorphism for every local
coefficient system $\mathcal B$ and every integer $n\ge d$.
\end{definition}

Our main result is

\begin{theorem}
Let $X$ be a connected CW--complex. The cohomology of $X$ is periodic
if and only if there is a spherical fibration $E\ra X$
with a total space $E$ that is homotopy equivalent to a finite
dimensional CW--complex.
\end{theorem}

In particular if $\Gamma$ is a {\sl discrete} group, then 
a $CW$--complex of type $K(\Gamma , 1)$ 
has periodic cohomology if and only if $\Gamma$ acts freely 
and properly on a
finite dimensional complex 
(i.e. a free $\Gamma$--CW complex)
that is homotopy equivalent to a sphere.
 From this we obtain a characterization of groups which act freely and
properly on $\mathbb R^n\times \mathbb \bbS^m$ which extends a conjecture
formulated by Wall \cite{W2} for groups of finite virtual
cohomological dimension:

\begin{corollary}
A discrete group $\Gamma$ acts freely and properly on
$\mathbb R^n\times\mathbb \bbS^m$ for some $m,n>0$ if and only if 
$\Gamma$ is a countable group with periodic cohomology.
\end{corollary}

Our results apply to natural situations arising from finite group
actions.  We consider actions that have isotropy subgroups with
periodic cohomology.

\begin{theorem}
Let $G$ be a finite group and let $X$ be a finite dimensional $G$--CW
complex such that every abelian subgroup of an isotropy subgroup is
cyclic.  Then for some integer $N>0$ there is a finite dimensional
CW--complex $Y$ such that $Y\simeq \bbS^N\times X$ and a free action
of $G$ on $Y$ such that the projection $Y\to X$ is $G$-equivariant. If
$X$ is simply--connected and
finitely dominated then $Y$ can be taken to be a finite
complex.
\end{theorem}

The rank $r(G)$ of a finite group $G$ is the dimension of its largest
elementary abelian subgroup. It has been conjectured (see \cite{BC1})
that Swan's theorem admits the following generalization: if $k$ is the
smallest integer such that a finite group $G$ will act freely on a
finite complex $Y\simeq\bbS^{n_1} \times\dots\times\bbS^{n_k}$, then
$k=r(G)$.  Using linear spheres, we prove:

\begin{theorem}
Let $G$ denote a finite $p$--group with center $Z(G)$ such that
$r(G)-r(Z(G))\le 1$; then $G$ acts freely on a finite complex
$Y\simeq \bbS^N\times (\bbS^{2[G:Z(G)]-1})^{r(G)-1}$.
\end{theorem}

As a corollary we obtain a succinct geometric 
characterization of low--rank $p$--groups which represents a long-sought
extension of Swan's theorem beyond spherical space forms:

\begin{corollary}
Let $G$ denote a finite $p$--group. Then
$G$ acts freely on a finite complex $Y\simeq \bbS^n\times\bbS^m$
if and only if $G$ does not contain a subgroup isomorphic to
$\bbZ/p\times\bbZ/p\times\bbZ/p$.
\end{corollary}

Constructing free actions of simple groups is a particularly difficult
problem. Using the classification of finite simple groups \cite{Go}
and a result due to Oliver \cite{O} we show there is no free action of
a rank two simple group on a product of spheres $\bbS^n\times\bbS^m$
that is a product action. Using character theory, we construct free
actions of some simple groups which are therefore necessarily exotic
(i.e. non--product actions).

\begin{theorem}
Each of the simple groups $A_5$, $SL_3(\bbF_2)$, $U_3(3)$, and $U_3(4)$ 
acts on linear spheres with rank one isotropy.  Hence each acts freely
and exotically on a finite complex that is homotopy equivalent to a
product of spheres $\bbS^n\times\bbS^m$.
\end{theorem}

Unfortunately many of the rank two simple groups do not have linear
spheres with the necessary properties. To deal with such cases, we
develop a method to produce $p$--local Euler classes from local
subgroups. We then assemble the actions on $p$--local spheres to
produce an action on a homotopy sphere with an appropriate integral
Euler class. Concentrating on the prime $2$ (indeed every non--abelian
simple group contains a copy of $\bbZ/2\times\bbZ/2$) we use the poset
of elementary abelian subgroups to produce suitable Euler classes
$2$--locally. Using additional arguments for odd primes and making a
case-by-case analysis we obtain

\begin{theorem}
Let $G$ denote a rank two simple group 
other than $PSL_3(\bbF_p)$,
$p$ an odd prime.  
Then $G$ acts freely on a finite complex
$Y\simeq\bbS^n\times\bbS^m$.
\end{theorem}

Our methods also imply
\begin{theorem}
Let $G$ denote a rank two finite group such that all of its Sylow subgroups
are either normal in $G$, abelian or generalized quaternion. Then $G$ acts
freely on a finite complex $Y\simeq\bbS^n\times\bbS^m$.
\end{theorem}

We are grateful to the CRM-Barcelona for its hospitality 
and to W.Dwyer, G.Mislin, B.Oliver, W.L\"uck, R.Solomon,
G.Seitz, S.Smith and O.Talelli for their helpful comments.
We also thank the referee for his useful remarks.
The results in this paper have been announced in
\cite{AS}. 

\section{Periodic complexes}

In this section we generalize a result of Swan \cite{Sw} concerning
finite groups with periodic cohomology. When coefficients do not
appear in a (co)homology group we assume that they are in $\bbZ$,
without any twisting. We first prove a number of lemmas required for
our main results.

\begin{definition}\label{def-periodic-cohomology}
A CW--complex $X$ has \emph{periodic} cohomology if there is a
cohomology class $\a\in H^*(X)$ with $|\a |> 0$ and an
integer $d\ge 0$ such that the map
\[
\a\cup-\,\colon H^n_{loc}(X,\mathcal B)\to H^{n+|\a|}_{loc}(X,\mathcal
B)
\]
induced by the cup product with $\a$ is an isomorphism for every local
coefficient system $\mathcal B$ and every integer $n\ge d$.
\end{definition}

A spherical fibration is a Serre fibration $E\to B$ for which the
fiber is homotopy equivalent to a sphere $\bbS^m$.  A spherical
fibration is orientable if the action of $\pi_1(B)$ on $H^m(\bbS^m)$ is
trivial and an orientation of the spherical fibration is a choice of
generator for $H^m(\bbS^m)$.  The Euler class of an oriented spherical
fibration is the cohomology class in $H^{m+1}(B)$ that is the transgression
of the generator for $H^m(\bbS^m)$.

We recall some basic constructions in homotopy theory. Note that we
will be working in the homotopy category throughout, so diagrams and
constructions which we use should be understood accordingly.  We
recall the notion of the fiber join of two fibrations.

\begin{definition}
Let $E_1\to X$ and $E_2\to X$ be two fibrations over the same
base space $X$, and let $E_1\times_XE_2$ be their fibered product. 
The {\sl fiber join} of these fibrations, denoted
by $E_1*_XE_2\to X$, is defined as
the homotopy pushout of the diagram of fibrations over $X$

\[
\begin{CD}
E_1\times_X E_2 @>p_2>> E_2\\
@Vp_1VV    \\
E_1  
\end{CD}
\]
where $p_1$ and $p_2$ are the natural projections.  The fibration
$E_1*_XE_2\to X$ has fiber $F_1*F_2$, the join of the respective
fibers of the two original fibrations.
\end{definition}

In what follows we will make use of basic properties of Postnikov sections 
(see \cite{Gray} for background).
The $q$--th Postnikov section of the $m$-sphere is denoted
$P_q\bbS^m$. Recall that the induced map $H^m(P_q\bbS^m) \to
H^m(\bbS^m)$ is an isomorphism if $q\ge m$.

\begin{definition}
Let $X$ be a connected CW--complex and let $k\ge0$ be an integer. A
cohomology class $\a\in H^n(X)$ is a \emph{$k$-partial Euler class} if
there is a fibration sequence
\[
P_{k+n-1}\bbS^{n-1}\to E \to X
\]
such that the action of $\pi_1(X)$ on $H^{n-1}(P_{k+n-1}\bbS^{n-1})=\bbZ$
is trivial and a generator of this group
transgresses
to $\a$.
\end{definition}

\begin{proposition}\label{prop-trivial-action}
For an integer $q\ge m$, let $P_q\bbS^m\to E\to B$ be a fibration
sequence with a connected base. If the action of $\pi_1(B)$ on
$H^m(P_q\bbS^m)$ is trivial then the action of $\pi_1(B)$ on all the
groups $H^*(P_q\bbS^m)$ is trivial.
\end{proposition}

\begin{proof}
For any fibration sequence $F\to E\to B$ there is a natural map
$\pi_1(B)\to [F,F]$ and the action of $\pi_1(B)$ on $H^*(F)$ is 
determined by the
composition $\pi_1(B)\to[F,F]\to Hom(H^*(F),H^*(F))$. The functor $P_m$
induces a map $[F,F]\to [P_mF,P_mF]$; in the proposition,
$F=P_q\bbS^m$, $P_mF=K(\bbZ,m)$ and $[P_mF,P_mF]=\bbZ$. So
the functor $P_m$ induces an isomorphism $[F,F]\to \bbZ$ 
(indeed, for $q\ge m$, we have isomorphisms $[\bbS^m,\bbS^m]
\cong [P_q\bbS^m,P_q\bbS^m]$ by the fundamental property
of Postnikov sections).
The action of $\pi_1(B)$ on $H^m(P_q\bbS^m)$ determines the
homomorphism $\pi_1(B)\to[F,F]=\bbZ$; if the action is trivial then the
homomorphism is trivial and the action of $\pi_1(B)$ on $H^*(P_q\bbS^m)$
is trivial.
\end{proof}

\begin{lemma}\label{lem-powers}
Let $X$ be a connected CW--complex and let $k\ge 0$ be an integer. For
every cohomology class $\a\in H^*(X)$ of positive degree there is an
integer $q\ge1$ such that the cup power $\a^q$ is a $k$--partial
Euler class.
\end{lemma}

\begin{proof}   
The lemma is proved by induction on $k$.  For $k=0$, let $f\colon X\to
K(\bbZ,n)$ be a map that represents the cohomology class $\a$ and let
$E\to X$ be the fibration defined by the pullback square
\[
\begin{CD}
E @>>> PK(\bbZ,n)\\
@VVV @VVpV\\
X @>f>> K(\bbZ,n)
\end{CD} 
\]
where $p$ is the path fibration.  The fiber of $E\to X$ is $K(\bbZ,n-1)$
which is $P_{n-1}\bbS^{n-1}$. By construction the
generator of $H^{n-1}(K(\bbZ,n-1))$ transgresses to $\a\in
H^*(X)$. Hence every $\a$ is a $0$--partial Euler class.

The inductive step is next. Assuming that $\a\in H^n(X)$ is a $k$--partial
Euler class, we show that some cup power $\a^q$ is a $k+1$--partial
Euler class.  Let
\[ 
P_{k+n-1}\bbS^{n-1}\to E\to X
\] 
be a fibration sequence for which a generator of
$H^{n-1}(P_{k+n-1}\bbS^{n-1})$ transgresses to $\a$.  By
Proposition~\ref{prop-trivial-action}, the action of $\pi_1(X)$ on
$H^*(P_{k+n-1}\bbS^{n-1})$ is trivial. So, the $E_2$ term of the Serre
spectral sequence for cohomology with coefficients in
$\pi_{k+n}(\bbS^{n-1})$ is
\[
E_2^{p,q}=H^p(X,H^q(P_{k+n-1}\bbS^{n-1},\pi_{n+k}(\bbS^{n-1}))).
\]
Let $g_k\colon P_{k+n-1}\bbS^{n-1}\to K(\pi_{k+n}(\bbS^{n-1}),k+n+1)$
be the $k$-invariant and let $\gamma_k$ denote the cohomology class
represented by $g_k$. The class $\gamma_k$ lies in $E_{2}^{0,k+n+1}$
and survives to $E_{k+3}$. The differential $d_{k+3}\gamma_k$ is a
class in $E_{k+3}^{k+3,n-1}$; if $k$ is sufficiently large, this group
is isomorphic to $H^{k+3}(X,\pi_{n+k}(\bbS^{n-1}))$ and $\a\cup
d_{k+3}\gamma_k=0$. If $d_{k+3}\gamma_k=0$ then $\gamma_k$ is
transgressive and $d_{k+n+2}\gamma_k\in H^*(X)/(\a)$.

Consider the diagram
\[
\begin{CD}
*^qP_{k+n-1}\bbS^{n-1} @>>> *^q_X E_k @>>> X\\
@VVV @VVV @VVV\\
P_{k+qn-1}(*^qP_{k+n-1}\bbS^{n-1}) @>>> P^q_XE_k @>>> X
\end{CD}
\]
having as its first row the $q$-fold fiberwise join of the fibration
$E_k \to X$ and as its second row the fibration obtained by the
fiberwise application of the Postnikov section functor $P_{k+qn-1}$ to
the first row. Applying the Postnikov section functor to the natural
map $\bbS^{qn-1}\to *^qP_{k+n-1}\bbS^{n-1}$ gives a weak equivalence
\[
P_{k+qn-1}\bbS^{qn-1}\to P_{k+qn-1}(*^qP_{k+n-1}\bbS^{n-1}).
\]
In dimension $k+qn+1$ one has 
\begin{eqnarray*}
H_{k+qn+1}(*^q P_{k+n-1}\bbS^{n-1})&=&\oplus_q \pi_{k+n}(\bbS^{n-1})\\
H_{k+qn+1}(P_{k+qn-1}\bbS^{qn-1})&=&\pi_{k+qn}(\bbS^{qn-1}).
\end{eqnarray*}

The induced map in homology is the $q$-fold sum of maps
$f_i: \pi_{k+n}(\bbS^{n-1})\to \pi_{k+qn}(\bbS^{qn-1})$.  By symmetry 
$f_i=f_j$ for all $1\le i,j \le q$.  To understand 
$f_1$ consider the composite
\[
P_{k+n-1}\bbS^{n-1}*(*^{q-1}\bbS^{n-1})\to *^q P_{k+n-1}\bbS^{n-1} 
\to P_{k+qn-1}\bbS^{qn-1}.
\]
The adjoint is the map 
\[
P_{k+n-1}\bbS^{n-1} \to \Omega^{(q-1)n}P_{k+qn-1}\bbS^{qn-1}=P_{k+n-1}
\Omega^{(q-1)n}\bbS^{qn-1}
\]
which is determined on the bottom cell. It is obtained by applying the
functor $P_{k+n-1}$ to the inclusion $\bbS^{n-1} \to
\Omega^{(q-1)n}\bbS^{qn-1}$. It follows that $f_1$ is the suspension
homomorphism and therefore that the map $\oplus f_i$ is the fold map
followed by the suspension homomorphism.

Let $g_k^{[q]}$ denote the $k$-invariant
\[
P_{k+qn-1}\bbS^{qn-1}\to K(\pi_{k+qn}(\bbS^{qn-1}),k+qn+1),
\] 
and let $\gamma_k^{[q]}$ be the cohomology class that it
represents. By naturality of differentials it follows
that $d_{k+3}\gamma_k^{[q]}$ is the image of $qd_{k+3}\gamma_k$ under the map
\[
H^{k+3}(X,\pi_{k+n}(\bbS^{n-1}))\to H^{k+3}(X,\pi_{k+qn}(\bbS^{qn-1}))
\]
that is induced by the suspension homomorphism on the
coefficients. 

Now the image of the suspension homomorphism is a finite group
and if $q$ is the exponent of the image then
$d_{k+3}\gamma_k^{[q]}=0$. On the other hand, if
$d_{k+3}\gamma_k=0$ then $\gamma_k$ transgresses to a class $\delta\in
H^{k+n+2}(X)$ which is well defined modulo the ideal generated by
$\a$ and $\gamma_k^{[q]}$ transgresses to $q \a^{q-1}
\overline{\delta}$ modulo the ideal generated by $\a^q$ where
$\overline{\delta}$ is the image of $\delta$ under the map
\[
H^{k+n+2}(X,\pi_{k+n}(\bbS^{n-1}))\to H^{k+n+2}(X,\pi_{k+qn}(\bbS^{qn-1}))
\]
induced by the suspension homomorphism.  Again, if $q$ is the exponent
of the image of the suspension homomorphism then $\gamma_k^{[q]}$
transgresses to zero modulo the ideal generated by $\a^q$. 

Let us now {\it choose} $q$ to be the square of the exponent of the image of the
suspension homomorphism. Then by construction,
$\gamma^{[q]}_k$ survives to $E_\infty$ in
the Serre spectral sequence. This implies that the map $g_k^{[q]}$ extends to
a map on $P_X^qE_k$ and there is a diagram
\[
\begin{CD}
P_{k+qn}\bbS^{qn-1} @>>> E_{k+1} @>>> X \\
@VVV @VVV @VVV\\
P_{k+qn-1}\bbS^{qn-1} @>>> P^q_XE_k@>>>  X \\
@VVV @VVV @VVV\\
K(\pi_{k+qn}(\bbS^{qn-1}),k+qn+1) @>>> K(\pi_{k+qn}(\bbS^{qn-1}),k+qn+1) @>>> *
\end{CD}
\]
in which each row and each column is a fibration sequence. Taking the 
fibration sequence in the first row completes the inductive step.
\end{proof}

\begin{lemma}\label{lem-spherical-fibration}
Let $P_{k+|\a|-1}\bbS^{|\a|-1}\to E \to X$ be a fibration sequence such
that a generator of $H^{|\a|-1}(P_{k+|\a|-1}\bbS^{|\a|-1})$ transgresses to
$\a\in H^*(X)$. If $k\ge d-2$ and the map
$\a\cup-\,\colon H^m(X)\to H^{m+|\a|}(X)$ is an isomorphism in integral
cohomology for $m\ge d$ then there is an orientable spherical
fibration $\overline E\to X$ with Euler class $\a$.  If in addition
the map $\a\cup-\,\colon H^m_{loc}(X,\mathcal B)\to
H^{m+|\a|}_{loc}(X,\mathcal B)$ is an isomorphism for every local
coefficient system $\mathcal B$ and every integer $m\ge d$ then
$\overline{E}$ is homotopy equivalent to a CW--complex of dimension
less than $d+|\a|$.
\end{lemma}

\begin{proof}
Let $n=|\a |$; then 
for each integer $t\ge k+n-1$ we will construct by induction a
fibration sequence $P_t\bbS^{n-1}\to E_t \to X$ and a map $E_{t+1}\to
E_t$ over $X$. The induction begins at $t=k+n-1$ with the fibration
sequence $P_{k+n-1}\bbS^{n-1}\to E \to X$ that one has by hypothesis.

The inductive step is next.  Suppose that one has a fibration sequence
$P_t\bbS^{n-1}\to E_t \to X$. The $k$-invariant $g_t\colon P_t\bbS^{n-1}\to
K(\pi_{t+1}(\bbS^{n-1}),t+2)$ extends to a map on $E_t$ if and only if the
cohomology class represented by $g_t$ survives
to $E_\infty$ in the Serre spectral sequence for cohomology with
coefficients in $\pi_{t+1}(\bbS^{n-1})$. There are two differentials
which might be non-trivial on the class represented by $g_t$. These are
$d_{t-n+4}$, taking values in the kernel of multiplication by $\a$; and
$d_{t+3}$, taking values in $H^*(X)/(\a)$. But by hypothesis both obstruction 
groups are zero and one can construct the diagram
\[
\begin{CD}
P_{t+1}\bbS^{n-1} @>>> P_t\bbS^{n-1} @>g_t>> K(\pi_{t+1}(\bbS^{n-1}),t+2)\\
@VVV          @VVV                @VVV\\
E_{t+1}       @>>>  E_t        @>h>> K(\pi_{t+1}(\bbS^{n-1}),t+2)\\
@VVV           @VVV            @VVV\\
X @>>> X      @>>>  *
\end{CD}
\]
in which each row and each column is a fibration sequence.

Let $E_\infty=\holim_t E_t$.  The homotopy fiber of the map
$E_\infty\to X$ is weakly equivalent to $\bbS^{n-1}$. By construction
the integral cohomology of $E_\infty$ vanishes above dimension $d+|\a
|-1$.  If in addition one has the condition on cohomology with local
coefficients that is given in the hypotheses, then it follows from the
Serre spectral sequence with local coefficients that $E_\infty$ has
finite cohomological dimension. So by a result of Wall \cite{W1},
$E_\infty$ is homotopy equivalent to a finite dimensional CW
complex. From this we obtain a spherical fibration over $X$ with total
space $\overline E$ homotopy equivalent to a complex of dimension less
than $d+|\a |$.
\end{proof}

\begin{lemma}\label{lem-killing-obstruction}
Let $X$ be a CW--complex and let $E\to X$ be an
orientable fibration such that the
fiber is equivalent to an odd dimensional sphere other than
$\bbS^1$, so $\pi_1(E)=\pi_1(X)$. For each integer $q\ge1$, the $q$-fold
fiber join $*^q_X E\to X$ is a spherical fibration.  If $E$ is
finitely dominated then $*_X^qE$ is finitely dominated and if
$\mathcal O_E\in\widetilde K_0(\bbZ\pi_1(X))$ is the finiteness obstruction of
$E$ then $q\mathcal O_E$ is the finiteness obstruction of $*^q_X E$.
\end{lemma}

\begin{proof}
Let $E_1\to X$ and $E_2\to X$ be orientable
spherical fibrations with fibers that
are odd dimensional spheres of dimension greater than one. Then there
is a pushout square
\[
\begin{CD}
E_1\times_X E_2 @>p_2>> E_2\\
@Vp_1VV    @VVV\\
E_1 @>>> E_1*_X E_2.
\end{CD}
\]
The projections $p_1$ and $p_2$ are spherical fibrations with odd
dimensional fibers.  If $E_1$ and $E_2$ are finitely dominated then
the fiber product $E_1\times_X
E_2$ and the fiber join $E_1*_X E_2$ are finitely dominated. It
follows that one has an equation in $\widetilde K_0(\bbZ\pi_1(X))$:
\[
\mathcal O_{E_1} + \mathcal
O_{E_2}=\mathcal O_{E_1*_X E_2}+\mathcal O_{E_1\times_X E_2}.
\]

Since $E_1\times_X E_2\to E_2$ is a spherical fibration with a finitely
dominated base and an odd dimensional sphere as the fiber, we conclude that
$\mathcal O_{E_1\times_X E_2}=0$ (see \cite{Lu}).  Repeated application 
allows us to prove the lemma.
\end{proof}

\begin{lemma}\label{lem-product-fibration}
Let $E \to X$ be a 
spherical fibration, where $X$ is a CW--complex
with k-skeleton $sk_kX$. For every integer $k\ge0$ there
is an integer $q$ such that the fibration $E_k\to sk_kX$ 
in the cartesian square
\[
\begin{CD}
E_k @>>> *^q_XE\\
@VVV @VVV\\
sk_kX @>>> X
\end{CD}
\]
is a product fibration.
\end{lemma}

\begin{proof}
First we observe that by taking fiber joins if necessary
(i.e. $E*_XE\to X$), we
can assume that $E\to X$ is an orientable
spherical fibration. 

Let $F$ be a CW--complex. Recall that $\aut F$ denotes the topological
monoid of self homotopy equivalences of $F$ and that $\aut_I F$ is the
connected component of the identity in $\aut F$. By \cite{BGM}, $B\aut_IF$
classifies fibrations $F\to E\to B$ for which the map
$\pi_1(B)\to [F,F]$ is trivial. In particular, $B\aut_I \bbS^m$ classifies
orientable spherical fibrations. 

The lemma will be proved by induction on $k$; the case
$k=0$ being trivial. 
Now, assume a trivialization over $sk_{k-1}X$; the obstructions to
extending this to a trivialization over $sk_kX$ are classes
${\mathcal O}_{\Delta}\in \pi_{k-1}(\aut_I(F))$, where $F\simeq
\mathbb S^m$ is the fiber and $\Delta$ ranges over the $k$--cells in
$X$.  Taking the $q$--fold fiber join replaces this obstruction with
$q\Sigma_*({\mathcal O}_{\Delta})$, where
$$\Sigma : Aut_I(F)\to Aut_I(*^q F)$$ is the natural stabilization map
defined by $f\mapsto f*1*\dots *1$.

Now we know that for a fixed $r>0$, $\pi_r(Aut_I(S^m))\cong
\pi^{st}_r(S^0)$ provided $m$ is sufficiently large. Note that these
stable homotopy groups of spheres are all finite. Hence given $k$ as
above, we can choose a large integer $q>0$ such that
$q\pi_{k-1}(Aut_I(*^qF))=0$, ensuring the vanishing of all the
obstructions and from there the extension of the trivialization.
\end{proof}

\begin{theorem}\label{thm-finite-dimensional-E}
The cohomology of a connected CW--complex $X$ is periodic if and only
if there is a spherical fibration $E \to X$ with a 
total space $E$ that is homotopy equivalent to a finite dimensional
CW--complex.
\end{theorem}

\begin{proof}
By taking fiber joins if necessary, 
we can assume that the spherical fibration is orientable.
Now if $E\to X$ is an oriented spherical fibration and the total space is
homotopy equivalent to a finite dimensional complex then the Gysin
sequence with local coefficients shows that $X$ has periodic
cohomology.

Conversely, assume that $X$ has periodic cohomology. Let $\a\in H^*(X)$
be an integral cohomology class and $d\ge0$ be an integer such that
$\a\cup-\,\colon H^m(X)\to H^{m+|\a|}(X)$ is an isomorphism for $m\ge
d$. By Lemma~\ref{lem-powers}, there is an integer $q\ge1$ and a fibration
sequence
\[
P_{d+q|\a|-1}\bbS^{q|\a|-1}\to E\to X
\]
such that a generator of $H^{q|\a|-1}(P_{d+q|\a|-1}\-\bbS^{q|\a|-1})$
transgresses to $\a^q\in H^{q|\a|}(X)$. By
Lemma~\ref{lem-spherical-fibration}, since $d\ge d-2$ there is an
orientable spherical fibration with Euler class $\a^q$, and the total
space has the homotopy type of a finite dimensional complex.
\end{proof}

Let $\Gamma$ denote a discrete group and let $X$ denote a
$CW$--complex which is an Eilenberg--MacLane space
of type $K(\Gamma , 1)$. We shall say that the group
$\Gamma$ has periodic
cohomology if $X$ has periodic cohomology.

\begin{corollary}
A discrete group $\Gamma$
has periodic cohomology if and only if 
$\Gamma$ acts 
freely and properly on a 
finite dimensional complex homotopy equivalent to a sphere.
\end{corollary}

If a countable discrete group $\Gamma$ acts freely and properly
on a finite dimensional complex $Y$ homotopy equivalent to a sphere
(i.e. $Y$ is a free $\Gamma$--CW complex), then it
acts freely and properly
on $\mathbb R^n\times\mathbb S^m$ for some
$m,n>0$. Indeed $Y/\Gamma$ has countable homotopy groups, hence is
homotopic to a countable complex which in turn is homotopic to an
open submanifold $V$ in some Euclidean space; applying the h--cobordism
theorem we can infer that for sufficiently large $q$ we have a 
diffeomorphism 
$\widetilde V\times\bbR^q\cong\bbR^n\times\bbS^m$ for some
$m,n>0$ (this appears in \cite{MT}, lemma 5.4; see also
\cite{CP}, page 139). Conversely the existence of 
a free and proper $\Gamma$--action
on $\bbR^n\times\bbS^m$ implies that $\Gamma$ is countable (as $\bbR^n\times
\bbS^m$ is a separable metric space) and that $\Gamma$
has periodic cohomology (via the Gysin sequence). 
Hence we obtain

\begin{corollary}
A discrete group $\Gamma$ acts freely and properly on $\mathbb R^n\times
\mathbb S^m$ for some $m,n>0$ if and only if $\Gamma$ is countable
and has periodic cohomology.
\end{corollary}

\begin{remark}
This result represents a proof of a
generalized version of a conjecture due to
Wall \cite{W2} for groups of finite virtual cohomological dimension
(that case was verified in \cite{CP}). More recently this has been extended
to a larger class of discrete groups (see \cite{MT}).  It is worth
noting that the authors in \cite{MT} use an abstract notion of
periodicity, namely the natural equivalence of functors between
$H^i(\Gamma, ~~~~)$ and $H^{i+n}(\Gamma ,~~~)$ for all $i\ge
d$. Moreover they speculate that this form of periodicity should be
equivalent to the existence of a free and proper
action on a finite dimensional
homotopy sphere. In view of the results presented here, this question
is equivalent to showing that abstract periodicity is always induced
by cup product with a cohomology class. This is an interesting open
problem in group cohomology.
\end{remark}

By a result of Wall \cite{W1}, a finitely dominated CW--complex $X$ is homotopy
equivalent to a finite CW--complex if and only if its finiteness
obstruction vanishes; the finiteness obstruction is an element in the
reduced projective class group $\widetilde K_0(\bbZ\pi_1(X))$.

\begin{theorem}\label{thm-finite-E}
Let $X$ be a CW--complex of finite type such that the reduced
projective class group $\widetilde K_0(\bbZ\pi_1(X))$ is a torsion group. The
cohomology of $X$ is periodic if and only if there is a spherical
fibration $E \to X$ such that the total space $E$ is homotopy
equivalent to a finite complex.
\end{theorem}

\begin{proof}
As before we can assume orientability
for our fibration by taking fiber joins
if necessary. 
Now if $E\to X$ is an orientable
spherical fibration and $E$ is homotopy
equivalent to a finite complex then the Gysin sequence shows that $X$
has periodic cohomology. 

Conversely, suppose $X$ has periodic cohomology. By
Theorem~\ref{thm-finite-dimensional-E} there is an orientable
spherical fibration
$E\to X$ such that $E$ is homotopy equivalent to a finite dimensional
CW--complex.  By forming the fiber join of the map $E\to X$ with itself
if necessary we may assume that the fiber of the map $E\to X$ is a
homotopy sphere of odd dimension $2n-1\ge3$ and hence that
$\pi_1(E)=\pi_1(X)$. Since $X$ has finite type, $E$ has finite type. Hence
$E$ is finitely dominated and its finiteness obstruction is an element
$\mathcal O_E\in\widetilde K_0(\bbZ\pi_1(X))$. By hypothesis $\mathcal O_E$
has finite order $q$.  The total space of the
$q$-fold fiber join $*_X^qE\to X$ is also finitely dominated.  By
Lemma~\ref{lem-killing-obstruction}, the finiteness obstruction of
$*_X^qE$ is $q\mathcal O_E=0$. Hence $*_X^qE$ is homotopy equivalent
to a finite complex.
\end{proof}

\begin{corollary}
Let $\Gamma$ denote a discrete group such that 
there exists a $CW$--complex of finite type which is an
Eilenberg--MacLane space of type $K(\Gamma , 1)$
and such that
$\widetilde K_0(\mathbb Z\Gamma)$ is torsion. Then $\Gamma$ has
periodic cohomology if and only if it acts freely 
and properly on a $CW$--complex
$Y$ homotopy equivalent to a sphere such that $Y/\Gamma$ is a 
finite complex.
\end{corollary}

\begin{remark}
If $G$ is a finite group, $\widetilde K_0(\bbZ G)$ is a finite abelian
group and of course we can find a $K(G, 1)$ which is a $CW$--complex
of finite type. 
\end{remark}

\section{Group Actions with Periodic Isotropy}

Let $G$ denote a finite group, we will be interested
in actions of $G$ on spaces such that the isotropy subgroups have
periodic cohomology. 
The key fact is given by

\begin{proposition}\label{prop-per-coh}
Let $X$ denote a finite dimensional $G$--CW complex, where $G$ is a finite
group.
Then the cohomology of $X\times_GEG$ is periodic
if and only if all the isotropy subgroups have periodic cohomology.
\end{proposition}
\begin{proof}
Assume that all the isotropy subgroups have periodic cohomology; then
we can choose $\a\in H^*(G,\bbZ)$ such that for any isotropy subgroup
$H\subset G$, $res^G_H(\a)$ is a periodicity generator (see \cite{A}, \cite{BC1}).
If $\a$ is described by the extension
$$0\to L_{\a}\to \Omega^r(\bbZ)\to \bbZ\to 0$$
where $\Omega^r(\bbZ)$ denotes the $r$--fold dimension shift of the
trivial $\bbZ G$--module $\bbZ$, then 
this is equivalent to the projectivity of the
modules $L_{\a}\otimes\bbZ [G/H]$ (see \cite{Be}, Chapter 5).

Now let $C^*(X)$ denote the cellular co-chains
on $X$, and consider the exact sequence
$$
0\to C^*(X)\otimes L_\alpha\to C^*(X)\otimes\Omega^{r}(\bbZ)\to C^*(X)\to
0.
$$
By construction, $C^*(X)\otimes L_\alpha$ is $\bbZ G$-projective, hence
cohomologically trivial---this yields a cohomology isomorphism
with any coefficients $M$ for all $i>dim~X$
$$
H^i(X\times_GEG,M){\buildrel\cong\over\longrightarrow} 
H^{i+|\a|}(X\times_GEG,M)
$$
which by construction is induced by multiplying with $\pi^*(\alpha
)$, where $\pi:X\times_GEG\to BG$ is the usual bundle map. 

For the converse observe that periodicity implies that for any
prime $p$, $H^*(X\times_GEG,\bbF_p)$ will have Krull dimension one
or zero,
hence by a theorem due to Quillen \cite{Q}, all the isotropy must
be of rank one i.e. with periodic cohomology.
\end{proof}

The preceding proposition
allows us to describe important specific situations where the
results from the previous section will apply. In particular,
the following theorem
is a
direct consequence of Theorem~\ref{thm-finite-dimensional-E}
and Lemma~\ref{lem-product-fibration}:

\begin{theorem}\label{thm-rank-one-isotropy}
Let $X$ denote a
finite dimensional $G$--CW complex ($G$ a finite group)
such that
all of its isotropy subgroups have periodic cohomology.
Then there
exists a finite dimensional CW complex $Y$ with a free
$G$--action such that $Y\simeq \bbS^N\times X$. If
$X$ is simply connected and finitely dominated, then
we can assume that $Y$ is a finite complex.
\end{theorem}

\begin{example}
Let $M_g$ denote a Riemann surface (orientable) of
genus $g>1$.  It is well known that $G=\hbox{Aut}(M_g)$ is a finite group (of
order at most $84(g-1)$) which acts preserving orientation on $M_g$ with
isolated singular points (see \cite{FK}). 
Furthermore it is known that the isotropy
subgroups of the action are cyclic.
Hence we obtain that $G=\hbox{Aut}~(M_g)$ acts freely on a
finite dimensional complex $Y\simeq \bbS^{2n-1}\times M_g$.
\end{example}

\begin{remark}
A more direct proof of Theorem~\ref{thm-rank-one-isotropy} can be given as follows.
Using induction on the dimension of the cells, one can construct a
spherical fibration over the skeleta of $X$ one stage at a time. 
The key ingredients
are the existence and homotopy uniqueness of the spherical space forms
for the isotropy subgroups and the coherence from the
choice of a global
periodicity generator; in fact the class $\pi^*(\alpha)$ constructed 
in Proposition~\ref{prop-per-coh} 
can be used. The rest of the proof consists of killing obstructions
(all of finite exponent)
by using fiber joins to extend through higher skeleta, and taking further
fiber joins to split the fibration (non-equivariantly). 
This approach is based on the method that 
was outlined in \cite{CP}.
\end{remark}

\begin{definition}
The {\sl $p$--rank} $r_p(G)$
of a finite group $G$ is the maximal
rank of an elementary abelian $p$--subgroup of $G$, and the
rank $r(G)$ of $G$ is $\max_{p||G|} r_p(G)$. 
\end{definition}

Let $G$ be a finite group of rank two. One can choose integral
cohomology classes $\a$, $\b$ in $H^*(BG)$ which form a homogeneous
system of parameters (see \cite{BC2}). This can be interpreted in turn 
as saying that
the module $L_{\a}\otimes L_\b$ is a projective $\bbZ G$-module. 
Under certain conditions this can be used to construct interesting
group actions.

\begin{theorem}\label{thm-rank-two-groups}
Let $G$ be a finite group of rank two. If $\a,\b\in H^*(BG,\bbZ)$ are
integral cohomology classes which form a homogeneous system of
parameters where $\b$ (of even degree larger than $2$) is the Euler
class of an oriented spherical fibration $X \to BG$, then there is a
free action of $G$ on a finite CW--complex $P$ such that $P$ is
homotopy equivalent to a product $\bbS^n\times \bbS^m$ and in the
Serre spectral sequence of the fibration $P\times_GEG\to BG$ the
generator of $H^n(\bbS^n,\bbZ)$ transgresses to $\b$ and the generator
of $H^m(\bbS^m,\bbZ)$ transgresses to a cup power of $\a$.
\end{theorem}

\begin{proof}
Let $X\to BG$ be an oriented spherical fibration having $\b$ as its
Euler class. Since $G$ is finite, $BG$ and hence $X$ are of finite
type. Moreover, 
since $L_\a\otimes L_\b$ is a projective $G$-module we infer from
the Gysin sequence that the map
$\a\cup-\colon H^s(X, M)\to H^{s+|\a|}(X, M)$ is an isomorphism for 
every coefficient $M$ and all integers $s$ which are larger than a
fixed integer $d$. By
Theorem~\ref{thm-finite-E} there is a spherical fibration $E\to X$
such that $E$ is homotopy equivalent to a finite CW--complex and the
fundamental class of the fiber transgresses to a cup power of $\a$. The
$G$ cover of $X$ is homotopy equivalent to a sphere $\bbS^n$ and so the
map $\widetilde X \to X$ factors up to homotopy through the
$n$-skeleton of $X$. By Lemma~\ref{lem-product-fibration} there is an
integer $q$ such that the spherical fibration $*^q_XE\to X$ is a
product fibration when restricted to the $n$-skeleton of $X$. Hence
the map $P\to \widetilde X$ in the cartesian square
\[
\begin{CD}
P @>>> *_X^qE\\
@VVV @VVV\\
\widetilde X @>>> X
\end{CD}
\]
is a product fibration and $P\to \widetilde X$ is a $G$-equivariant
spherical fibration such that the fundamental class of the fiber
transgresses to a cup power of $\a$. So the complex $P$ has the
required properties.
\end{proof}

 From the above we infer that to construct a free action of a rank two
finite group on a homotopy product of spheres all we need is an
Euler class $\b$ which forms part of a homogeneous system of parameters.
We deal with this question in the next section.

\section{Euler Classes from Actions on Spheres}

In this section we will investigate the existence of group actions on
homotopy spheres which have Euler classes representing chosen
cohomology classes. Here we assume that all groups which appear are
finite.
Let $\beta\in H^N(G)$ denote a cohomology class
which happens to be the Euler class of a $G$-action on $X$, a CW
complex homotopic to the $(N-1)$-sphere.  In this case the algebraic
and topological Gysin sequences provide an identification
$$
H^i(X\times_GEG,\bbZ)\cong H^{i+1}(G,L_\beta)\qquad \forall \; i\ge N.
$$

We recall (see \cite{Be} for background)

\begin{definition}
Let $M$ denote any finitely generated $\bbZ G$--module
which is $\bbZ$--torsion free (a lattice).
The {\sl complexity} $cx_G(M)$ of $M$ is the maximal rate of growth
of $H^*(E, M\otimes\bbF_p)$  
taken over all $p$--elementary abelian subgroups of $G$ and
primes $p$ dividing $|G|$. 
\end{definition}

If $\beta$ is an Euler class, then
the complexity of $L_{\b}$ is determined by the asymptotic
growth rates of the $H^*(E,L_{\beta}\otimes \bbF_p)$, and so we have

\begin{proposition}
If $\b$ is an Euler class for a $G$--CW complex
$X\simeq \bbS^{N-1}$, then 
$cx_G(L_{\b})$ is precisely the maximum of the Krull dimensions of
$H^*(X\times_GEG,\bbF_p)$, as $p$ ranges over all prime divisors
of $|G|$.
\end{proposition}

\begin{definition} A cohomology class
$\beta\in H^N(G,\bbZ)$ is said to be {\it
$p$--effective\/} if $cx_G(L_{\b}\otimes\bbF_p)<r(G)$. 
The class $\b$ is said to be {\it effective} if it is
$p$--effective for all primes $p$ dividing $|G|$ or equivalently
if $cx_G(L_{\b})<r(G)$.
\end{definition}

\begin{remark}
Note that for a class to be effective it suffices to check that it
is $p$--effective for primes of maximal rank.
\end{remark}

We assume from now on that $\b$ is an effective Euler class.  Using
the lemma above, this will have geometric implications if $X$ is {\sl
finite dimensional}, an assumption we shall make for the rest of this
section.  It implies
that the maximal rank elementary abelian subgroups
in $G$ cannot  
act\footnote{Of course it is well-known that if $(\bbZ /p)^r$ acts on a
space such as $X$, then there is an index $p$ subgroup fixing a point
in $X$.}with stationary points.
In fact we have a simple characterization of effective Euler
classes:

\begin{lemma}\label{lemma-effective}  
If a cohomology class $\beta\in H^N(G)$ is the
Euler class of an action on a finite dimensional
$X\simeq \bbS^{N-1}$, then $\beta$ is
effective if and only if every maximal rank elementary abelian
subgroup of $G$ acts without stationary points.
\end{lemma}

Using the localization theorem in \cite{Q}, it is clear
that this condition can be detected homologically, i.e.
$\beta$ is effective if and only if 
$\beta\big|_E\not= 0$ $\forall \; E\subset G$ elementary
abelian of maximal rank.

\begin{example}
Let $E=(\bbZ /2)^n$ and $V$ the complex reduced
regular representation of $E$, which is $(2^n-1)$-dimensional.  Then
$V^E=\{ 0 \}$ and we obtain an Euler class $e$ associated to the action on
$S(V)\cong \bbS^{2(2^n-1)-1}$.  In this case $e\in H^{2(2^n-1)}(E)$ is
$D^2$, the square of the top Dickson class $D=\prod_{0\not= z\in H^1(E,\bbF_2)}z$.
Note that any proper subgroup of $E$ fixes a point in $S(V)$, i.e.\
$e\big|_{E'}=0$ $\forall\; E'\;\subsetneq E$.
\end{example}

Euler classes arising from representations are of course fundamental and
will play a key role in some of our results.  
Indeed for $p$--groups it turns out that representations
can always be used:

\begin{theorem} 
Let $G$ denote a finite $p$--group acting on a
finite complex $X$ homotopic to a sphere. Then there exists a real representation
$V$ of $G$ such that the Euler class associated to $X$ is effective if and
only if the Euler class for $S(V)$ is effective.
\end{theorem}

\begin{proof} 
We apply the main result in \cite{DH}, namely given a complex as
above, there exists a real representation $V$ of $G$ such that 
$dim~V^H=dim~X^H+1$ for all subgroups $H\subset G$. Applying 
Lemma~\ref{lemma-effective} completes
the proof.
\end{proof}

 From our characterization of effective Euler classes, we see
that the problem of finding an effective `linear' Euler class 
reduces to finding a representation $V$ such that $V^E=\{ 0\}$
for all $E\subset G$ maximal elementary abelian $p$-subgroups. 
We provide a construction of such representations for certain groups.

\begin{proposition}\label{center-rep} 
Let $G$ denote a finite group with center 
$Z(G)\cong \bbZ/n_1\times\dots\times\bbZ/n_k$ and
corresponding generators $x_1,\dots ,x_k$.
Then there exist $[G:Z(G)]$--dimensional complex
representations
$V_1, V_2,\dots ,V_k$ of $G$ such that $V_j^H=\{ 0\}$
for any subgroup $H\subset G$ containing $x_j$.
\end{proposition}

\begin{proof}
Let $\chi_j :Z(G)\to\bbC$ denote the 
one-dimensional representation of $Z(G)$ defined
by $x_j\mapsto e^{2\pi i/n_j}$, and $x_h\mapsto 1$ for $h\ne j$.
Let $V_j=\hbox{Ind}^G_{Z(G)}(\chi_j)$.  Using the double coset
formula and the fact that $x_j$ is central, we readily see that
$V_j\big|_{Z(G)}\cong\bigoplus^{[G:Z(G)]}\chi_j$ and $V_j^{\langle
x_j\rangle}=\{ 0\}$.
\end{proof}

If $X=S(V_1)\times \cdots \times S(V_k)$ with
the diagonal $G$-action, then $Z(G)$ acts freely
on $X$.  In fact it follows from the results in \cite{D} that if the
largest central elementary abelian $p$--subgroup of $G$ has rank equal
to $r$, then the associated
Euler classes $e_1,\ldots ,e_r\in H^{2[G:Z(G)]}(G,\bbF_p)$ 
form a regular sequence in cohomology.

We use these representations\footnote{Note that if 
$G$ is a $p$--group such that every element of order $p$
is central, the above construction yields a free $G$--action on
$X\cong (\bbS^{2[G:Z(G)]-1})^{r(G)}$.} to show

\begin{proposition}\label{center-one}
Let $G$ denote a finite group such that
the $p$--rank of its center
$Z(G)$ is equal to one less than the rank of $r(G)$ of $G$. 
Then $G$ acts on $X=(\bbS^{2[G:Z(G)]-1})^{r(G)-1}$
such that the isotropy subgroups are all of $p$--rank at most equal
to one. 
\end{proposition}

\begin{proof} 
Take $X$ as before, constructed from
the representations induced from
the center. If any $E\cong\bbZ/p\times \bbZ/p\subset G$ fixed a point in
$X$, then $E\cap Z(G)=\{1\}$. Therefore the subgroup generated by
$E$ and $Z(G)$ would have $p$--rank at least as large as
$r(G)+1$, a contradiction.
\end{proof}

Recalling that a $p$--group always has non--trivial center, we obtain

\begin{corollary}\label{rank-two-sphere}  
If $G$ is a finite $p$-group of rank
equal to two, then $G$ acts on $X=\bbS^{2[G:Z(G)]-1}$ with rank one isotropy
subgroups.
\end{corollary}

Applying Theorem~\ref{thm-rank-one-isotropy} to the
$G$--CW complex $X$ (see \cite{tD}), we obtain

\begin{theorem} 
Let $G$ denote a $p$-group such that $r=$
rank $G\le$ rank $Z(G)+1$.  Then $G$ acts freely on a finite complex
$$
Y\simeq \bbS^N\times (\bbS^{2[G:Z(G)]-1})^{r-1}.
$$
\end{theorem}

In particular, combining the well--known fact \cite{H}
that
$(\bbZ/p)^3$ does not act freely on $\bbS^n\times\bbS^m$
with Corollary~\ref{rank-two-sphere}, 
we obtain a geometric chacterization of rank two $p$--groups extending
Swan's result for periodic groups

\begin{theorem}
A finite $p$--group $G$ acts freely on a finite complex
$Y\simeq\bbS^n\times \bbS^m$
if and only if it does not contain a subgroup isomorphic to 
$\bbZ/p\times\bbZ/p\times\bbZ/p$.
\end{theorem}

In the next section we will use linear Euler classes to construct
actions of simple groups. It is natural to ask if such spheres will
suffice for any rank two group of composite order. Unfortunately this
is not the case, in fact the symmetric group $\Sigma_5$ has no
effective linear Euler classes. One could try instead to use
non--linear actions as was done in the spherical space form
problem\footnote{For example the non--abelian group of order $21$ does
not admit a fixed--point free representation; as shown by Petrie in
\cite{P} there is a smooth non--linear action on the $5$--dimensional
sphere with Euler class a periodicity generator.}, but this approach
is not systematic enough to yield complete results.

\section{Examples: Linear Euler Classes and Actions of Simple Groups.}

In this section we will construct actions of a number of rank two groups
on a product of  two spheres.  The initial input is a representation
sphere $S(V)$ which has an Euler class which is {\it effective}.  We will
concentrate on groups which are of rank two at the prime $p=2$, in
particular on simple groups. 

To put these results in context, we recall a few facts about simple
groups.  
First, note that by the Feit-Thompson Theorem (\cite{Go}, pg.15), 
we know that non-abelian
simple groups are of even order.  Less well-known
is that fact that any non-abelian simple group must contain $\bbZ /2\times
\bbZ /2$ as a subgroup (see \cite{Go}, pg.13)\footnote{This of course explains
why simple groups do not arise in studying spherical space forms.}.
Furthermore there is an explicit classification of
simple groups $G$ such that their rank at $p=2$ is exactly two 
(see \cite{Go}, pg.72). Examining this list and using well--known facts
about finite groups of Lie type (see \cite{Car}) we obtain the following
{\sl complete} list of rank two simple groups i.e. simple groups which
do not contain $(\bbZ/p)^3$ for {\sl any} prime number $p$:
$$
\hbox{PSL}_2(\bbF_p), p\;\hbox{ odd }\; p\ge 5;\; 
\hbox{PSL}_2(\bbF_{p^2}), p\;\hbox{ odd };\; \hbox{PSL}_3(\bbF_p),\;
p\hbox{ odd };$$
$$\hbox{PSU}_3(\bbF_p),\; p\hbox{ odd }; \;\hbox{PSU}_3(\bbF_4);\; A_7\;
\hbox{and}~~~M_{11}.
$$

On the other hand, one can verify\footnote{We are grateful to R. Solomon
for explaining this to us.} that all the groups listed above
contain a copy of $A_4$, the alternating group on four letters.
This can be used to obtain a very basic result on actions of
rank two simple groups:

\begin{theorem} 
If $G$ is a simple group of rank equal to two,
then it cannot act freely on a finite dimensional complex
$X\simeq \bbS^{n_1}\times \bbS^{n_2}$ via a product
action.
\end{theorem}

\begin{proof}
Suppose such an action exists. Recall that
the $k$--fold join of $\bbS^r$ is $\bbS^{k(r+1)-1}$.
Hence taking joins we can obtain
a free action of $A_4$
on $X'=\bbS^{(n_1+1)(n_2+1)-1}\times \bbS^{(n_1+1)(n_2+1)-1}$, which contradicts
a result due to Oliver \cite{O}.
\end{proof}

Hence to construct actions of these groups one must necessarily consider
exotic, non--product actions that go beyond products of representations.
It is known (see \cite{R}) that no simple group
can act freely on a product of spheres arising from representations,
but of course Oliver's criterion and its natural generalization
to the groups $(\bbZ/2)^n\times_T\bbZ/(2^n-1)$
(see \cite{Si}) seem more effective in a topological setting.

Consider a group $G$ which contains $\bbZ /2\times \bbZ/2$ as a subgroup. We
are interested in constructing actions of $G$ on spheres such that the
Euler class is effective in mod 2 cohomology. 
The complex reduced regular representation of $\bbZ /2\times \bbZ/2$
gives rise to an orientation--preserving
action on $\bbS^5$ with effective mod 2 Euler class.
To extend this we
make use of the following elementary lemma. 

\begin{lemma}
Let $\chi$ be a complex character for $G$ which
takes
the constant value $\chi (2)$ on all 
involutions.
Then $\chi\big|_E$ is
the character associated to a multiple of the
reduced regular representation of $E\cong
\bbZ /2\times \bbZ /2$, for every
$E\subseteq G$ if and only if $\chi (1)=-3\chi (2)$.
\end{lemma}

\begin{proof}
The irreducible representations of $E\cong \bbZ /2\times \bbZ
/2$ are described as follows:  if $E=\langle a,b\rangle$ let
$$
\chi_1(a)=\chi_1(b)=1, ~~~~~
\chi_2(a)=-1,\quad \chi_2(b)=1$$
$$\chi_3(a)=1,\quad \chi_3(b)=-1,~~~~
\chi_4(a)=-1,\quad \chi_4(b)=-1.$$

Then $\chi\big|_E=m_1\chi_1+m_2\chi_2+m_3\chi_3+m_4\chi_4$.  Evaluating on
$a, b$, $ab$ respectively yields
$$
m_1-m_2+m_3-m_4 =
m_1+m_2-m_3-m_4 =
m_1-m_2-m_3+m_4 =\chi (2)$$
whence we readily infer that $m_2=m_3=m_4$ and $m_1=\chi (2)+m_2$.  Hence
we deduce that $m_1=0$ if and only if $m_2=-\chi (2)$; but noting that
$\chi (1)=m_1+3m_2$ we conclude that $m_1=0$ if and only if $\chi
(1)=-3\chi (2)$.

\end{proof}

Inspecting the character tables in the Atlas \cite{Co}, we see that the
groups $\Sigma_4$, $A_5$ and $\hbox{SL}_3(\bbF_2)$ have $3$-dimensional
complex characters which satisfy the conditions of the lemma above.

\begin{theorem}  
$\Sigma_4$, $A_5$ and $\hbox{SL}_3(\bbF_2)$ act
freely on finite complexes $Y\simeq \bbS^N\times \bbS^5$. 
\end{theorem}

\begin{proof}  
Note that $|A_5|=2^2\times 5\times 3$, $|\Sigma_4|=2^3\times
3$, $|\hbox{SL}_3(\bbF_2)|=2^3\times 3\times 7$, hence for odd
primes $p$ they have cyclic $p$-Sylow subgroup and we need only
consider the prime $p=2$.  By our lemma and
the existence of the characters
mentioned above, the free actions on $Y\simeq \bbS^N\times \bbS^5$ can be
constructed using Theorem~\ref{thm-rank-one-isotropy}. 
\end{proof}

Next we consider slightly more complicated simple groups.

\begin{theorem}
The simple group $U_3(3)$ acts freely on a
finite complex $Y_1\simeq\bbS^{N_1}\times \bbS^{11}$; whereas the simple
group $U_3(4)$ acts freely on $Y_2\simeq \bbS^{N_2}\times \bbS^{23}$.
\end{theorem}

\begin{proof}
We first consider $G=U_3(3)$, note that it has order
$2^5\cdot 3^3\cdot 7$ and that $\hbox{rank}_2(G)=\hbox{rank}_3(G)=2$.  Now
$P=\hbox{Syl}_3(P)$ is an extra-special $3$-group 
\[
1\to \bbZ /3\to P\to \bbZ
/3\times \bbZ /3\to 1
\]
with $\bbZ /3\subset P$ central.  Denote by $3A$ the
conjugacy class of this central element of order 3.

 From the Atlas \cite{Co} we see that $G$ has one conjugacy class of involutions
$(2A)$ and two conjugacy classes of elements of order three ($(3A)$ and
$(3B)$ respectively). There is a complex character $\chi_2$ which takes
the following values:
\[
\table{}
! 1 ! 2A ! 3A ! 3B\rr
$\chi_2$ ! 6 ! -2 ! -3 ! 0
\caption{}
\]
 From these values and the previous lemma, we infer that the associated mod
2 Euler class is effective.  However notice that the subgroup generated by
3A has no fixed-points, hence acts freely on the associated sphere
$S(\chi_2)\cong \bbS^{11}$.  Therefore in mod 3 cohomology the Euler class
restricts non-trivially to the cohomology of 3A and so must be
effective.  We deduce that $G$ acts on $\bbS^{11}$ with rank one isotropy and
so by our theorem $G$ acts freely on a finite complex $Y_1\simeq
\bbS^{N_1}\times \bbS^{11}$.

Now we consider $G=U_3(4)$; in this case $|G|=2^6\cdot 3\cdot 5^2\cdot
13$.  According to the Atlas \cite{Co}, $G$ has only one conjugacy class of
involutions, 2A, and furthermore there is a complex character $\chi$ with
the following values
$$\table{}
! 1 ! 2A ! 5A\rr
$\chi$ ! 12 ! -4 ! -3
\caption{}
$$
This implies that the mod 2 Euler class for $S(\chi )\cong \bbS^{23}$ is
effective, and as before (but now for $p=5$) we also see that 5A acts
without fixed-points, hence the Euler class is mod 5 effective too.  Hence
$G$ acts on $\bbS^{23}$ with periodic isotropy and the rest follows as
before.
\end{proof}

\begin{example} 
$G=\hbox{GL}_2(\bbF_3)$, $|G|=2^4\cdot 3$.  We
consider a portion of its character table, with a complete set of
conjugacy classes:
$$\table{}
! 1 ! 2A ! 4 ! 2B ! 3\rr
$\chi_2$ ! 1 ! 1 ! 1 ! -1! 1\rr
$\chi_6$ ! 2 ! -2 ! 0 ! 0 ! -1\rr
$\chi_8$ ! 4 ! -4 ! 0 ! 0 ! 1 
\caption{}
$$
It is not hard to see that $\hbox{GL}_2(3)$ acts freely on
$S(\chi_6)\times S(\chi_2)\cong \bbS^3\times \bbS^1$.  From \cite{AM} we can
think of $H^*(\hbox{GL}_2(3), \bbF_2)$ as the subalgebra of $\bbF_2[x_1, y_1]$
generated by $d_1=x_1+y_1$, $d_4=(x_1y_1)^2$, $d_3=x_1^2y_1+x_1y_1^2$,
namely
$$
\bbF_2 [d_1, d_4, d_3]\big/ d^2_3=d_4d^2_1.
$$
$S(\chi_6)$ has Euler class equal to $d_4$, whereas $S(\chi_2)$ has Euler
class equal to $d^2_1$.

We can, however, construct an exotic action as follows.  Let
$$
\chi =2\chi_2+\chi_8
$$
then $\chi (2A)=\chi (2B)=-2$, $\chi (1)=6$.  Applying our lemma, we obtain
an action of $\hbox{GL}_2(3)$ on $\bbS^{11}$ with periodic isotropy and
hence a free action on $Y=\bbS^N\times \bbS^{11}$.  In this case the class
associated to the $\bbS^{11}$ is $d^4_3$, and as universal periodicity
class $\alpha$ we can take $d^4_1+d_4$, hence $N=4\ell -1$ and we get
$\alpha^\ell =(d^4_1+d_4)^\ell$.
This exotic action will be used in the sequel.
\end{example}

\section{Local Euler Classes and Local Subgroups}
 
As we have seen, linear spheres will often have effective Euler
classes.  However, there are many groups where this is not the case,
the first interesting example is $\Sigma_5$, the symmetric group on
five letters.  Similarly most of the rank 2 simple groups previously
listed must be handled in a different way.

To deal with this we introduce local Euler classes and a notion of local
equivalence between groups which is flexible enough to substantially
extend our collection of examples, at least locally.
In our next section we show that effective local Euler classes can be
assembled to provide integral effective Euler classes.

Let $\bbS^m_{(p)}$ denote the $p$-local $m$-sphere; so $\bbS^1_{(p)}$
is the classifying space $B\bbZ_{(p)}$ of the $p$-local integers and for
$m>1$, $\bbS^m_{(p)}$ is a simply connected space with
$H_m(\bbS^m_{(p)},\bbZ)=\bbZ_{(p)}$.

\begin{definition}
Let $X$ be a connected CW--complex and let $k\ge0$ be an integer. A
cohomology class $\a\in H^n(X, \mathbb Z_{(p)})$ 
is a \emph{$p$--local Euler class}
if there is a fibration sequence
$\bbS^{n-1}_{(p)}\to E \to X$
such that the action of $\pi_1(X)$ on 
$H^{n-1}(\bbS^{n-1}_{(p)}, \mathbb Z_{(p)})=\bbZ_{(p)}$
is trivial and a generator of 
$H^{n-1}(\bbS^{n-1}_{(p)}, \mathbb Z_{(p)})$ transgresses
to $\a$.
\end{definition}

\begin{definition}
Let $G$ and $H$ denote discrete groups. We say that they are
$p$--equivalent (denoted $G\simeq_p H$) for a prime $p$ if the
localizations at $p$ 
of their classifying spaces
are homotopy equivalent, $BG_{(p)}\simeq
BH_{(p)}$. 
\end{definition}

Note that if we are given two $p$--equivalent groups such that one of
them is finite, then the notion of $p$--effective Euler class makes
sense for both of them.

\begin{definition}
A $p$--effective Euler class for a discrete group $G$ is a 
cohomology element $x\in H^{2n}(BG,\mathbb Z_{(p)})$ which
is of the form $f^*(y)$, where $f:BG_{(p)}\to BH_{(p)}$ is
an $H\mathbb Z_{(p)}$--equivalence and 
$y\in H^{2n}(BH_{(p)}, \mathbb Z_{(p)})$ is a $p$--effective
Euler class for the finite group $H$.
\end{definition}

We will be interested in constructing {\sl local} actions, i.e. actions on
the localization $\bbS^{N-1}_{(p)}$ with effective mod $p$ Euler classes.
A basic result is explained in the following lemma, which will
be proved in the following section.
 
\begin{lemma}\label{p-equiv}
Assume that we have two discrete groups $G$ and $H$ which are
$p$--equivalent, where one of them is finite.  Then there exists an
action of $G$ on $\bbS^{2k-1}_{(p)}$ (a $p$--local sphere) with a
$p$--effective Euler class if and only if the same holds for $H$.
\end{lemma}

\begin{example}
Using a double coset decomposition one can see that the natural
inclusion $\Sigma_{p^n}\subset\Sigma_{p^n+1}$ ($p$ a prime) induces an
isomorphism in mod $p$ cohomology.  Hence $\Sigma_5$ will act on a
$2$--local sphere with effective Euler class (as this holds for
$\Sigma_4$).
\end{example}

\begin{example}
If $G$ is a finite group with abelian $p$--Sylow subgroup $A$, then
(by a theorem of Swan), the inclusion of the normalizer $N_G(A)$ in
$G$ induces a mod $p$ cohomology equivalence.  (see \cite{AM}).
\end{example}

\begin{proposition}\label{syl-ab}
If a finite group $G$ has a normal $p$--Sylow subgroup, then its
cohomology has a $p$--effective Euler class. Hence if $G$ has an
abelian $p$--Sylow subgroup it has a $p$--effective Euler class.
\end{proposition}

\begin{proof}
Let $G$ denote a group with normal $p$--Sylow subgroup 
$P\triangleleft G$.
Consider an effective class $\a$ for $P$ arising from a
$[P:Z(P)]$--dimensional representation of $P$ as in 
Proposition~\ref{center-rep}. Now induce this to a representation for
$G$; by the double coset formula, the resulting Euler class will
restrict in $H^*(P,\bbF_p)$ to $\prod_{g\in G/P}~g\alpha$, which is evidently also an
effective Euler class.
\end{proof}

Next we describe a method for producing characters for discrete groups
which arise as the fundamental group $\Gamma$ of a graph of finite
groups (see \cite{Br}, p. 179 and \cite{Se}). Recall that associated
to such a group we have a finite connected graph $Y$ together with
groups $G_v$ and $G_e$ where $v$ (respectively $e$) ranges over the
vertices (respectively edges) of $Y$. We also have monomorphisms
$\theta_0:G_e\to G_v$, $\theta_1:G_e\to G_w$ for every edge $e$ with
vertices $v$ and $w$. There is a tree $T$ with an action of $\Gamma$
such that the orbit space of the action is precisely $Y$. The group
$\Gamma$ can be thought of as a generalized amalgam obtained as a
quotient of the free product of the vertex groups $G_v$ by the
smallest normal subgroup generated by the relations
$\theta_0(a)=\theta_1(a)$ for all $a\in G_e$.

Now fix a group $Q$, and suppose we have homomorphisms
$f_v:G_v\to Q$ such that $f_v\theta_0=f_w\theta_1$ for all
pairs of vertices inciding on a common edge. This uniquely determines
a homomorphism $F:\Gamma\to Q$. In particular if $Q=GL_n(\bbC)$, we
can use this approach to construct linear representations for
generalized amalgams.

The following theorem makes this useful for constructing Euler
classes for a rank two group $G$. Let $|A_p(G)|$ denote the geometric
realization of the poset of non--trivial $p$--elementary abelian
subgroups in $G$ (see \cite{AM} for properties). Note that conjugation
induces a natural $G$--CW complex structure on this space.

\begin{proposition}\label{poset-space}
Let $G$ denote a finite group such that $|A_p(G)|$ is a
path--connected graph.  If $\Gamma$ denotes the fundamental group of
the graph of groups associated to $Y=|A_p(G)|/G$, then the natural
surjection $p: \Gamma\to G$ induces a $\bbZ_{(p)}$--cohomology
isomorphism.
\end{proposition}

This is a special case of a theorem due to Brown (see \cite{Br} and
also \cite{We}).  The key fact is that $H^1(|A_p(G)|,\bbZ_{(p)})$ is a
projective $\bbZ_{(p)}G$--module (analogous to the Steinberg module).
It is worthwhile to point out that $\pi=ker~p=\pi_1(|A_p(G)|)$ is a
free group of finite rank and that $\Gamma$ can be identified with
$\pi_1(|A_p(G)|\times_GEG)$. In this case the tree $T$ is simply the
universal cover $\widetilde{|A_p(G)|}$. Note that $|A_p(G)|$ is path
connected if and only if $\Gamma\to G$ is surjective. 

More generally,
(as explained in \cite{We}, pg. 151-152)
$G$ acts transitively on the components of
$|A_p(G)|$, with a stabilizer $N$ which is
$p$--equivalent to $G$ via the restriction map. Hence for our
purposes we can
always replace $G$ with this stabilizer group and so apply
the previous result
to construct Euler classes for a rank two group by using a 
$p$--equivalent group $\Gamma$. In other words we make use of the
equivalences $G\simeq_p N\simeq_p\Gamma$.

Our strategy will be to glue together characters on local subgroups to
yield characters for generalized amalgams, and then use the associated
localized spheres to produce effective Euler classes.  Note that the
Euler class associated to a $\Gamma$--sphere $W$ will be effective if
and only if it restricts to an effective Euler class on each of the
generating finite subgroups. Hence we need to assemble characters for
the generalized amalgams which arise from characters with effective
$\bbZ_{(p)}$--Euler classes.

\begin{example}
$G=M_{11}$, $p=2$.  In this case (see \cite{AM}) we can express
$M_{11}$ as a quotient of $\Gamma =\hbox{GL}_2(3)*_{D_8} \Sigma_4$
with $\Gamma\twoheadrightarrow G$ inducing a mod 2 equivalence.
Consider $\chi =2\chi_2+\chi_8$ the representation for
$\hbox{GL}_2(3)$ which we used previously, and compare it with the
following character for $\Sigma_4$ (also used previously):
$$\table{}
! 1 ! (12) ! (12)(34) ! (1234) ! (123)\rr
$\nu$ ! 3 ! -1! -1! 1 ! 0
\caption{}
$$
Then clearly $\chi$ agrees with $2\nu$ on elements of order $2$ and $4$,
hence on the subgroup $D_8$.  Therefore
we have a representation $V$
($6$-dimensional) for $\Gamma$ which affords an effective Euler class mod
2 (indeed it restricts to effective classes on the two conjugacy classes
of maximal
finite subgroups). We obtain a $2$--effective Euler
class in $H^{12}(M_{11},\bbF_2)$. 
Note that 
$V=\bbZ /3\times \bbZ /3=Syl_3(M_{11})$, with normalizer of order $144$.
Hence we may find a $3$--local sphere of dimension $31$ providing a $3$--effective
Euler class in $H^{32}(M_{11},\bbF_3)$.
\end{example}

\begin{example}
$G=\hbox{PSL}_2(\bbF_q)$ $q$ odd, $q\equiv \pm 3$ mod
$8$.  In this case the $2$--Sylow subgroup is elementary abelian
and by Proposition~\ref{syl-ab} we can obtain a $2$--effective
Euler class.
For the case $G=\hbox{PSL}_2(\bbF_q)$ $q$ odd, $q\equiv \pm 1$ mod
$8$, we observe that the group has two conjugacy classes of
rank two elementary abelian subgroups (both
with normalizers isomorphic to $\Sigma_4$), and one conjugacy class of
involutions, with centralizer a dihedral group $D_{2r}$ (here $r$ depends
on the prime $p$). The corresponding intersections are both $D_8$, hence
we can easily construct a $3$--dimensional representation for all of these
stabilizers with the desired property. We obtain a $2$--effective Euler
class in $H^6(\hbox{PSL}_2(\bbF_q),\bbF_2)$.
\end{example}

\begin{example}
$G=\hbox{PSL}_3(\bbF_q)$ $q$ odd, $q\equiv -1$  mod
$3$.  Here the $2$--local structure is
$$\xymatrix{
&\hbox{GL}_2(q)~~{\buildrel{\bbZ /q-1\wr\bbZ/2}
\over{\bullet\!\!\hbox{\fib}\!\!\bullet}}~~(\bbZ/q-1\times\bbZ/q-1)\times_T\Sigma_3\\}
$$
Take $\bbZ /2\subset\bbZ /q-1\times \bbZ /q-1\subset\hbox{GL}_2(q)$ and choose
$\chi$ a $1$-dimensional representation of $\bbZ /q-1\times \bbZ /q-1$ such that
$\chi\big|_{\bbZ /2}$ is non-trivial.  Let $\nu=\hbox{Ind}^{GL_2(q)}_{\bbZ
/q-1\times \bbZ /q-1}(\chi )$; then this is a
$q(q+1)$-dimensional representation on which
$\hbox{SL}_2(q)$ acts $2$-freely.  Now taking
$\hbox{GL}_2(q)/\hbox{SL}_2(q)\cong \bbZ /q-1$ and selecting a
$1$-dimensional embedding $\bbZ /q-1\subset\bbS^1$ we obtain a character
$\omega$
such that every element of even order in $\hbox{GL}_2(q)$ acts freely on
$S(\nu)\times S(\omega)$.  Now take
$$
{\mathcal K} =2\nu+q(q+1)\omega.
$$
This is a $3q(q+1)$-dimensional representation.  Let 2A denote the central
involution, 2B the non-central one, then
$$
{\mathcal K} (2A)=-q(q+1), \,\,
{\mathcal K} (2B)=-q(q+1)
$$
and hence $S({\mathcal K})$ has an effective Euler class mod 2.

Next we consider $\hbox{SL}_3(\bbF_q)\supset (\bbZ /q-1\times \bbZ
/q-1)\times_T\Sigma_3$, where this action arises from $\bbZ
/q-1\wr\Sigma_3\subset\hbox{GL}_3(\bbF_q)$ using the fact that diagonal
matrices of determinant one are permuted by $\Sigma_3$.

Now let $\bbZ /q-1{\buildrel\rho\over\hookrightarrow} U(1)$ be the standard
$1$-dimensional representation given by taking a generator to a primitive
$(q-1)$-root of unity; this gives rise to a representation
$$
(\bbZ /q-1)^2\times_T\Sigma_3\hookrightarrow U(1)\wr \Sigma_3\subset U(3)
$$
and so we obtain a $3$-dimensional character $\zeta$ for $(\bbZ /q-1\times
\bbZ /q-1)\times_T\Sigma_3$.  By naturality $S(\zeta )$ has an effective
$6$-dimensional Euler class which restricts to the square of the top
Dickson class on all $\bbZ /2\times \bbZ /2$ subgroups --- indeed one checks
that $\zeta\big|_{\bbZ /2\times \bbZ /2}$ is a copy of the reduced regular
representation.  It is direct to verify that $q(q+1)\zeta$
agrees with ${\mathcal K}$
on $\bbZ/q-l\wr\bbZ/2$, 
and so we can construct a $6q(q+1)$-1 dimensional
sphere associated to the amalgam $\Gamma$ providing an effective mod 2
cohomology class. Observe that if $q\equiv 1$ mod $3$, then suitably modified
arguments produce a similar result, where we must divide out by the
central $\bbZ/3$ in the diagram above.
\end{example}

\begin{example}
For the groups $U_3(q)$, $q\equiv 1$ mod $3$,
we have the following diagram of subgroups\footnote{Again we are 
grateful to R.Solomon
for sketching a description of $A_2(U_3(p))$.} 

$$\xymatrix{
&\hbox{GU}_2(q)~~{\buildrel{\bbZ /q+1\wr\bbZ/2}
\over{\bullet\!\!\hbox{\fib}\!\!\bullet}}~~(\bbZ/q+1\times\bbZ/q+1)
\times_T\Sigma_3\\}
$$
A construction similar to the previous one produces the desired
linear sphere for the amalgam. In the case $q\equiv -1$ mod $3$ we
must divide by a central $\bbZ/3$ in the diagram above and modify the
construction accordingly.
\end{example}

The inclusion $A_6\hookrightarrow A_7$ induces a mod $2$ equivalence,
hence given the isomorphism $A_6\cong PSL_2(\bbF_9)$ we see that $A_7$
has an effective $2$--local Euler class.  Thus we have completed
arguments for all the simple rank 2 groups on our list, and our
results can be collected together as follows:

\begin{theorem}
Let $G$ denote a finite simple group of rank equal to two. Then $G$
acts on some $Y\simeq \bbS^{2k-1}_{(2)}$ such that the associated Euler
class is $2$--effective.
\end{theorem}

Odd primes are in most instances easier to handle than $p=2$. If $q$
is an odd prime, then the $q$--Sylow subgroups of the $PSL_2(\bbF_p)$
and $PSL_2(\bbF_{p^2})$ are all cyclic (as the corresponding $SL_2$
have rank one). In addition the odd order Sylow subgroups for $A_7$
and $M_{11}$ are all abelian and we have already seen how to construct
locally effective Euler classes in that case.  The only complications
may arise either at $p=3$ or at the defining characteristic for
$U_3(p)$ and $L_3(p)$. It is however elementary to check that for any
rank $2$ simple group other than $U_3(3)$ or $L_3(3)$ the $3$--Sylow
subgroup is abelian--and of course we have already dealt with $U_3(3)$
in the previous section.  The following lemma takes care of $U_3(p)$
at the characteristic

\begin{lemma}
Let $G=U_3(p)$, the unitary group over $\bbF_p$. Then $G$ has a unique
cuspidal unipotent representation and it has degree $p(p-1)$.
Furthermore the associated Euler class in $H^{2p(p-1)}(G,\bbZ)$
is $p$--effective.
\end{lemma}

\begin{proof}
The first part can be found in \cite{L}, page 174. For the
second part we use an argument due to G.Seitz (\cite{Sei}).
Note that $p(p-1)$ is the minimal degree of a complex irreducible
character for $G$. This is also the minimal degree for a faithful
irreducible of the Sylow normalizer, which is a semidirect
product $P\times_T\bbZ/p-1$, where $P=Syl_p(G)$ is extra--special
of order $p^3$. The cyclic group $\bbZ/p-1$
acts transitively on the non--trivial linear characters of $Z(P)$.
The given representation restricts to $Syl_p(G)$ as the sum of the $p-1$
irreducibles of degree $p$. Restricting further to $Z(P)$ we just get 
all non--trivial linear characters, each $p$ times. This has value
$-p$ on each element of order $p$ in the center, in particular
any such element acts without fixed--points on the associated linear
sphere. We infer from Lemma~\ref{lemma-effective} that the corresponding
Euler class must be $p$--effective.
\end{proof}

To summarize the results in this section, we can state:

\begin{theorem}\label{simple-effective}
Given any rank two simple group $G$ and any prime $q$ dividing $|G|$
(except possibly $PSL_3(\bbF_p)$ at $p$), 
 then it has a
 $q$--effective Euler class in its cohomology.  \end{theorem}

\section{Homotopy Actions}

In this section we prove Lemma~\ref{p-equiv} and show that effective
Euler classes can be obtained by  gluing together $p$--effective Euler
classes.

Let $F$ be a CW--complex. Then $\aut F$ is the topological monoid of
self homotopy equivalences of $F$ and $\aut_I F$ is the component of
the identity in $\aut F$. By \cite{BGM} there is a universal fibration
with base the classifying space $B\aut F$ and homotopy fiber $F$ such
that, for $X$ a CW--complex, the fiber homotopy equivalence classes of
fibrations with homotopy fiber $F$ and base $X$ are in natural one to
one correspondence with the homotopy classes of maps from $X$ to
$B\aut F$. Let $G$ be a topological group, the equivariant homotopy
classes of free $G$-CW complexes are in a natural one-to-one
correspondence with the fiber homotopy equivalence classes of
fibrations with base $BG$. Combining these results, the homotopy
equivalence classes of maps $BG\to B\aut F$ are in one-to-one
correspondence with the homotopy actions of $G$ on $F$. 

For a set of primes $K$, let $\bbZ_K$ be the smallest subring of
the rationals that contains the reciprocals of the prime numbers not in
$K$. A simply-connected CW--complex is $\bbZ_K$-local if one of
two equivalent conditions hold: its integral homology groups are
$\bbZ_K$-modules or its homotopy groups are $\bbZ_K$-modules
(the equivalence comes from Serre $\mathcal C$--theory). The localization
of a simply-connected CW--complex $X$ is a CW--complex $X_K$ and a
map $X\to X_K$ (called the localization map) such that the $X_K$ is
$\bbZ_K$-local and the induced map $H_*(X,\bbZ)\otimes \bbZ_K \to
H_*(X_K,\bbZ)$ is an isomorphism.  

We recall a standard result in homotopy theory (see \cite{Dwyer}, 2.1). 
\begin{proposition}
Let $K$ denote any set of primes in $\mathbb Z$, then $B\aut_I \mathbb S^n_K$
is a $\bbZ_K$-local space.
\end{proposition}

Using the above we can now prove Lemma~\ref{p-equiv}.  Suppose given
two $p$--equivalent groups $G$ and $H$ (with one of them finite) and
an action of $G$ on $\bbS^n_{(p)}$ with $p$--effective Euler class. As
$B\aut_I \bbS^n_{(p)}$ is $p$--local, the corresponding map $BG\to
B\aut_I\bbS^n_{(p)}$ factors through $BG_{(p)}$. Now using the
$p$--equivalence and the natural localization map we obtain a homotopy
class in $[BH,B\aut_I(\bbS^n_{(p)})]$ which by naturality must give
rise to an action with $p$--effective Euler class.

Next we make a small digression to
outline how local homotopy actions can be assembled to define
a global action.
Let $J$ and $K$ be sets of prime
numbers. Then (see \cite{BK}) there is a homotopy pullback square of
localization maps
\[
\begin{CD}
X_{J\cup K} @>>> X_K \\
@VVV @VVV \\
X_J @>>> X_{J\cap K}.
\end{CD}
\]
The map $\hom(X_J,X_J)\to\hom(X,X_J)$ induced by the localization map is
a weak equivalence and so there is a homotopy pullback square
\[
\begin{CD} 
\hom(X_{J\cup K},X_{J\cup K}) @>>> \hom(X_K,X_K)\\
@VVV @VVV \\
\hom(X_J,X_J) @>>> \hom(X_{J\cap K},X_{J\cap K}).
\end{CD}
\]
The localization $f_{J\cup K}$ of a map $f$
is a weak equivalence if and only if the maps $f_J$ and
$f_K$ are weak equivalences and so there is a homotopy pullback square
of monoids
\[
\begin{CD} 
\aut X_{J\cup K} @>>> \aut X_K\\
@VVV @VVV \\
\aut X_J @>>> \aut X_{J\cap K}.
\end{CD}
\]
If $X$ is a finite nilpotent complex, the localization $f_{J\cup K}$
of the self map $f\colon X\to X$ is homotopic to the identity map on
$X_{J\cup K}$ if and only if the self maps $f_J$ and $f_K$ are homotopy
equivalent to the identity maps on $X_J$ and $X_K$. So, there is a
homotopy pullback square of monoids
\[
\begin{CD} 
\aut_I X_{J\cup K} @>>> \aut_I X_K\\
@VVV @VVV \\
\aut_I X_J @>>> \aut_I X_{J\cap K}.
\end{CD}
\]
Constructing the classifying spaces of the monoids gives a homotopy
pullback square
\[
\begin{CD} 
B\aut_I X_{J\cup K} @>>> B\aut_I X_K\\
@VVV @VVV \\
B\aut_I X_J @>>> B\aut_I X_{J\cap K}.
\end{CD}
\]
This pull-back construction 
will be used to assemble together local homotopy actions.

We are now ready to 
consider the problem of gluing together $p$--effective Euler
classes for a finite group $G$.
An integral cohomology class $\a\in H^n(BG,\bbZ)$
is an Euler class if there is a spherical fibration $\bbS^{n-1}\to E\to
BG$ such that the fundamental class in 
$H^{n-1}(\bbS^{n-1}, \mathbb Z)$ transgresses to
$\a$. It is an effective Euler class if its restriction to every
elementary abelian subgroup of maximal rank is non-trivial.  A
$p$-local cohomology class $\a\in H^n(BG,\bbZ_{(p)})$ is a a $p$-local
Euler class if there is a fibration $\bbS^{n-1}_{(p)}\to E\to BG$ such
that the fundamental class in 
$H^{n-1}(\bbS^{n-1}_{(p)},\bbZ_{(p)})$ transgresses
to $\a$. It is an effective $p$--local Euler class if its restriction to
every elementary abelian $p$-subgroup of maximal
rank is nontrivial.
In this section we show that the problem of constructing effective
Euler classes can be solved one prime at a time.

\begin{theorem}
An integral cohomology class $\a\in H^{2n}(BG,\bbZ)$ with $n\ge2$ is an
effective Euler class if and only if for each prime $p$ such that
$rkG=rk_pG$, its $p$-localization $\a_p\in
H^{2n}(BG,\bbZ_{(p)})$, is an effective $p$-local Euler class.
\end{theorem}

\begin{proof}
Clearly an Euler class is effective if and only if its
$p$-localization at each prime for which $rkG=rk_pG$ is
effective. Moreover if $\a$ is an Euler class then each of its
localizations $\a_p$ is a $p$-local Euler class.

Now assume that $\a\in H^{2m}(BG,\bbZ)$ and that for each prime $p$ for which
$rkG=rk_pG$ one has a fibration $E_p\to BG$ with homotopy fiber the
$p$-local sphere $\bbS^{2m-1}_{(p)}$ and $p$-local Euler class $\a_p$. So
$\pi_1(E_p)=G$ and the universal cover of $E_p$ is a $G$-CW complex that
is homotopy equivalent to $\bbS^{2m-1}_{(p)}$. Since the fundamental
class survives to $E_2$ in the Serre spectral sequence, the action of
$G$ on the $p$-local homotopy sphere is homologically trivial and
hence homotopically trivial.  Let $p_1,p_2,\dots,p_k$ be
the finite set of primes for which $rkG=rk_pG$ and let $K$ be the set
of all other primes. There is a unique homologically trivial action of
$G$ on the rationalization $\bbS^{2m-1}_0$ of an odd dimensional sphere
(here all primes are inverted).  Give
$\bbS^{2m-1}_K$ the trivial action of $G$. Then one has a diagram of
$G$-CW complexes
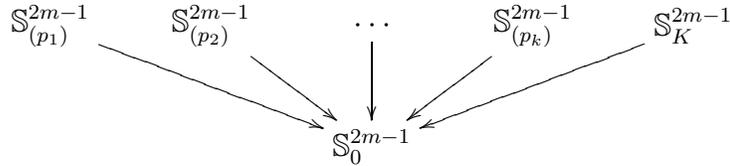
\[
\xymatrix{
\bbS^{2m-1}_{(p_1)}\ar[drr] & \bbS^{2m-1}_{(p_2)}\ar[dr] & {\cdots}\ar[d] & 
\bbS^{2m-1}_{(p_k)}\ar[dl] & \bbS^{2m-1}_K\ar[dll]\\
 & & \bbS^{2m-1}_0
}
\]

Using the homotopy pullback diagram for the classifying
spaces described previously, we see that the homotopy classes defining the
local homotopy actions define a homotopy action of $G$ on the
homotopy pull-back of the diagram, which by the
fracture lemma of Bousfield-Kan (see \cite{BK}, page 147)
is homotopy equivalent to
$\bbS^{2m-1}$. From this we obtain a $G$--CW complex
whose associated Borel construction gives the required spherical
fibration 
\[
\bbS^{2m-1}\to \bbS^{2m-1}\times_G EG\to BG.
\]
\end{proof}

Our methods yield actions on homotopy spheres which are hard to realize
geometrically. In particular they allow us to construct effective Euler
classes for groups where linear spheres are not useful. In turn we may use
these spheres to construct free group actions on finite complexes.
Using fiber joins and applying Theorem~\ref{syl-ab}
and Theorem~\ref{simple-effective} respectively, we obtain

\begin{theorem}
Let $G$ denote a rank two finite
group such that all of its $p$--Sylow subgroups
are either normal in $G$, abelian or generalized quaternion. 
Then $G$ acts freely on a finite complex $Y\simeq \bbS^n\times\bbS^m$
\end{theorem}

\begin{theorem}
Let $G$ denote a rank two finite simple group other than
$PSL_3(\bbF_p)$, $p$ an odd prime. Then $G$ acts freely on a finite complex
$Y\simeq \bbS^n\times\bbS^m$.
\end{theorem}

We conclude with an example illustrating the theorem above.
\begin{example}
Let $G=M_{11}$, the first Mathieu group. As we have seen it has a effective 
Euler classes for $p=2$ (in dimension $12$) and for $p=3$ (in dimension $32$),
where $p=2,3$ are the only two maximal rank primes. Hence we obtain
that $M_{11}$ acts freely on a finite complex $Y\simeq \bbS^N\times\bbS^{383}$
(note that $*^{32}\mathbb S^{11}\simeq *^{12}\mathbb S^{31}\simeq 
\mathbb S^{383}$).
\end{example}

\end{document}